\documentclass[10pt]{article}
\usepackage{amssymb,amsmath,amsthm,mathrsfs}

\newcommand{\ie}{\emph{i.e.}}

\newcommand{\cf}{\emph{cf}}

\newcommand{\Real}{\mathbb{R}}

\newcommand{\sii}{L^2}

\newcommand{\Sobi}{W_0^{1,2}}
\newcommand{\sobi}{W^{1,2}}
\newcommand{\Smooth}{C}
\newcommand{\Hilbert}{\mathcal{H}}
\newcommand{\Dom}{D}
\newcommand{\spec}{\mathrm{spec}\;\!}
\newcommand{\eps}{\varepsilon}
\newtheorem{Lemma}{Lemma}
\newtheorem{Theorem}{Theorem}
\newtheorem{Corollary}{Corollary}
\theoremstyle{remark}
\newtheorem*{Remark}{Remark}
\begin{document}
%
\title{\textbf{\Large
Hardy inequalities in strips on ruled surfaces
\footnote{
Accepted for publication in      
\emph{Journal of Inequalities and Applications}.
}
}}
\author{\textsc{David Krej\v{c}i\v{r}\'{\i}k}}
\date{\small 17 August 2005}
\maketitle
%
\begin{abstract}
\noindent
We consider the Dirichlet Laplacian
in infinite two-dimensional strips
defined as uniform tubular neighbourhoods
of curves on ruled surfaces.
We show that the negative Gauss curvature
of the ambient surface gives rise
to a Hardy inequality
and use this to prove certain stability
of spectrum in the case of
asymptotically straight strips
about mildly perturbed geodesics.
\end{abstract}
%

\section{Introduction}
%
Problems linking the geometry of two-dimensional manifolds
and the spectrum of associated Laplacians have been considered
for more than a century.
While classical motivations come from theories
of elasticity and electromagnetism,
the same rather simple models can be also remarkably successful
in describing even rather complicated phenomena
in quantum heterostructures.
Here an enormous amount of recent research has been undertaken
on both the theoretical and experimental aspects
of binding in curved strip-like waveguide systems.

More specifically, as a result of theoretical studies
it is well known now that
the Dirichlet Laplacian in an infinite planar strip
of uniform width always possesses
eigenvalues below its essential spectrum
whenever the strip is curved and asymptotically straight.
We refer to \cite{ES,GJ} for initial proofs
and to~\cite{DE,LCM,KKriz} for reviews
with many references on the topic.
The existence of the curvature-induced bound states
is interesting from several respects.
First of all, one deals with a purely quantum effect
of geometrical origin, with negative consequences
for the electronic transport in nanostructures.
From the mathematical point of view,
the strips represent a class
of non-compact non-complete manifolds
for which the spectral results of this type
are non-trivial, too.

At the same time,
a couple of results showing that
the attractive interaction due to bending
can be eliminated by appropriate additional perturbations
have been established quite recently.
Dittrich and K\v{r}{\'\i}\v{z} \cite{DKriz2}
demonstrated that the discrete spectrum of the Laplacian
in any asymptotically straight planar strip
is empty provided the curvature
of the boundary curves does not change sign
and the Dirichlet condition
on the locally shorter boundary
is replaced by the Neumann one.
A different proof of this result
and extension to Robin boundary conditions
were performed in~\cite{FK3}.
Ekholm and Kova\v{r}{\'\i}k~\cite{MK-Kov} obtained
the same conclusion for the purely Dirichlet Laplacian
in a mildly curved strip
by introducing a local magnetic field
perpendicular to the strip.
The purpose of the present paper is to show
that the same types of repulsive interaction
can be created if the ambient space of the strip
is a negatively curved manifold instead
of the Euclidean plane.

A spectral analysis of the Dirichlet Laplacian
in infinite strips embedded in curved
two-dimensional manifolds
was performed for the first time
by the present author in~\cite{K1}.
He derived a sufficient condition which guarantees
the existence of discrete eigenvalues
in asymptotically straight strips;
in particular, the bound states exist
in strips on positively curved surfaces
and in curved strips on flat surfaces.
He also performed heuristic considerations
suggesting that the discrete spectrum might be empty
for certain strips on negatively curved surfaces.
Similar conjectures were also made previously
for strips on ruled surfaces in~\cite{CB1}.
However, a rigorous treatment of the problem remained open.

In the present paper, we derive several Hardy inequalities
for mildly curved strips on ruled surfaces,
which proves the conjecture for this class of strips.
A ruled surface is generated by straight lines
translating along a curve in the Euclidean space;
hence its Gauss curvature is always non-positive.
The reason why we restrict to ruled surfaces
in this paper is due to the fact that
the Jacobi equation determining the metric
in geodesic coordinates is explicitly solvable,
so that rather simple formulae are available.
Nevertheless, it should be possible to extend
the present ideas to other classes of non-positively
curved surfaces for which more precise information
about geodesics are available.

Hardy inequalities represent a powerful technical tool
in more advanced theoretical studies of elliptic operators.
We refer to the book~\cite{Opic-Kufner}
for an exhaustive study and generalizations
of the original inequality due to Hardy.
Interesting Hardy inequalities
on non-compact Riemannian manifolds
were established in~\cite{Carron_1997}.
In the quantum-waveguide context,
various types of Hardy inequality
were derived in~\cite{MK-Kov,B-MK-Kov,EKK}
in order to prove certain stability of spectrum
of the Laplacian in tubular domains.

Here the last reference is the closest
to the issue of the present paper.
Indeed, the authors of~\cite{EKK} considered
a three-dimensional tube constructed
by translating a non-circular
two-dimensional cross-section
along an infinite curve
and obtained that the twisting
due to an appropriate construction
eliminates the curvature-induced discrete spectrum
in the regime of mild curvature.
Formally, the strips of the present paper can be viewed
as a singular case of~\cite{EKK}
when the cross-section is replaced by a segment
and the effect of twisting
is hidden in the curvature of the ambient space.
While~\cite{EKK} and the present paper
exhibit these similarity features,
and also the technical handling of the problems is similar,
they differ in some respects.
On the one hand, the present situation is simpler,
since it happens that the negative curvature
of the ambient space gives rise to an explicit
repulsive potential (\cf~(\ref{potential}) below)
which leads to a Hardy inequality
in a more direct way than in~\cite{EKK}.
On the other hand, we do not perform
the unitary transformation of~\cite{EKK}
in order to replace the Laplacian on the Hilbert space
of a curved strip by a Schr\"odinger-type operator
on a ``straighten'' Hilbert space,
but work directly with ``curved'' Hilbert spaces.
This technically more complicated approach
has an advantage that we need to impose
no conditions whatsoever on the derivatives of curvatures.

Although we are not aware of a direct physical
interpretation of the Laplacian in infinite strips
if the ambient space has a non-trivial curvature,
there exists an indirect motivation
coming from the theory of quantum layers
studied in~\cite{DEK2,CEK,LL1}.
In these references, the Dirichlet Laplacian
in tubular neighbourhoods of a surface
in the Euclidean space is used for the quantum Hamiltonian
(\cf~\cite{Encinosa-Mott_2003} for a similar model).
Taking our strip as the reference surface,
the layer model of course differs from the present one,
but a detailed study of the latter
is important to understand
certain spectral properties of the former.
Similar layer problems are also considered
in other areas of physics away from quantum theories,
\cf~\cite{Gridin-Craster-Adamou_2005}.
Finally, the present problem is
a mathematically interesting one
in the context of spectral geometry.

The organization of the paper is as follows.
The ambient ruled surface, the strip
and the corresponding Dirichlet Laplacian
are properly defined in the preliminary Section~\ref{Sec.Pre}.
In Section~\ref{Sec.geodesic}, we consider
the special situation of the strip
being straight in a generalized sense.
If the Gauss curvature of such a strip does not vanish
identically and the strip is thin enough,
we derive a central Hardy inequality
of the present paper, \cf~Theorem~\ref{Thm.Hardy}.
In fact, the latter is established by means
of a ``local'' Hardy inequality, \cf~(\ref{local}),
which might be also interesting for applications.
In Section~\ref{Sec.curved},
we apply Theorem~\ref{Thm.Hardy} to mildly curved strips
and prove certain stability of spectrum,
\cf~Theorem~\ref{Thm.stability}.
As an intermediate result,
we obtain a general Hardy inequality
for mildly curved strips on ruled surfaces,
\cf~(\ref{Hardy.c}).

\section{Preliminaries}\label{Sec.Pre}
%
Given two bounded continuous functions~$\kappa$ and~$\tau$
defined on~$\Real$ with~$\kappa$ being positive,
let $\Gamma:\Real\to\Real^3$
be the unit-speed curve whose curvature and torsion
are~$\kappa$ and~$\tau$, respectively.
$\Gamma$~is determined uniquely up to congruent transformations
and possesses a distinguished $\Smooth^1$-smooth Frenet frame
$\{\dot{\Gamma},N,B\}$ consisting of tangent,
normal and binormal vector fields,
respectively (\cf~\cite[Chap.~1]{Kli}).
It is also convenient to include
the case of~$\kappa$ and~$\tau$
being equal to zero identically,
which corresponds to~$\Gamma$
being a straight line with a constant Frenet frame.

Given a bounded $\Smooth^1$-smooth function~$\theta$ defined on~$\Real$,
let us introduce the mapping $\mathcal{L}:\Real^2\to\Real^3$ via
\begin{equation}\label{surface}
  \mathcal{L}(s,t)
  := \Gamma(s)
  + t\,\big[N(s)\;\!\cos\theta(s)-B(s)\;\!\sin\theta(s)\big]
  \,.
\end{equation}
$\mathcal{L}$~represents a ruled surface
(\cf~\cite[Def.~3.7.4]{Kli})
provided it is an immersion.
The latter is ensured by requiring that the metric tensor
$G\equiv(G_{ij})$ induced by~$\mathcal{L}$, \ie
$$
  G_{ij}
  := (\partial_i\mathcal{L})\cdot(\partial_j\mathcal{L})
  \,, \qquad
  i,j \in \{1,2\}
  \,,
$$
where the dot denotes the scalar product in~$\Real^3$,
is positive definite.
Employing the Serret-Frenet formulae (\cf~\cite[Sec.~1.3]{Kli}),
we find
\begin{equation}\label{metric}
  G =
  \begin{pmatrix}
    h^2 & 0 \\
    0   & 1
  \end{pmatrix}
  \,, \qquad
  h(s,t)
  := \sqrt{\big[1-t\;\!\kappa(s)\cos\theta(s)\big]^2
  + t^2 \;\! \big[\tau(s)-\dot{\theta}(s)\big]^2 }
  \,.
\end{equation}
Hence, it is enough to assume that~$t$ is sufficiently small
so that the first term in the square root defining~$h$
never vanishes.

More restrictively,
given a positive number~$a$,
we always assume
\begin{equation}\label{Ass.basic}
  a \, \|\kappa\cos\theta\|_\infty < 1
  \,,
\end{equation}
so that also~$h^{-1}$ is bounded,
and define a \emph{ruled strip}
of width~$2a$ to be the Riemannian manifold
\begin{equation}\label{strip}
  \Omega := \big(\Real \times (-a,a),G\big)
  \,.
\end{equation}
That is, $\Omega$~is a non-compact and non-complete surface
which is fully characterized by
the functions~$\kappa$, $\tau$, $\theta$
and the number~$a$.
It is easy to verify that the Gauss curvature~$K$ of~$\Omega$
is non-positive, namely,
\begin{equation}\label{Gauss}
  K =
  -(\tau-\dot{\theta})^2 \, h^{-4}
  \,.
\end{equation}
Moreover, if the mapping~$\mathcal{L}$ is injective,
then the image
$
  \mathcal{L}\big(\Real\times(-a,a)\big)
$
has indeed the geometrical meaning of
a non-self-intersecting strip
and~$\Omega$ represents its parameterization
in geodesic coordinates.

\begin{Remark}
In~(\ref{metric}),
let us write~$k$ instead of $\kappa\cos\theta$
and~$\sigma$ instead of $\tau-\dot{\theta}$,
and assume that~$k$ and~$\sigma$
are given bounded continuous functions on~$\Real$.
Then, abandoning the geometrical interpretation
in terms of ruled surfaces based on~$\Gamma$,
(\ref{strip})~can be considered
as an abstract Riemannian manifold,
with $a\|k\|_\infty < 1$ being the only restriction.
The spectral results of this paper extend automatically
to this more general situation
by applying the above identification.
\end{Remark}

Our object of interest is the Dirichlet Laplacian in~$\Omega$,
\ie, the unique self-adjoint operator~$-\Delta_D^\Omega$
associated with the closure of the quadratic form~$Q$
defined in the Hilbert space
\begin{equation}\label{Hilbert}
  \Hilbert := \sii(\Omega)
  \equiv \sii\big(\Real\times(-a,a),h(s,t)\,ds\,dt\big)
\end{equation}
by the prescription
\begin{equation}\label{form}
  Q[\psi]
  := \big(\partial_i\psi,G^{ij}\partial_j\psi\big)_\Hilbert
  \,, \qquad
  \psi\in\Dom(Q)
  := \Smooth_0^{\infty}\big(\Real\times(-a,a)\big)
  \,,
\end{equation}
where $(G^{ij}):=G^{-1}$
and the summation is assumed
over the indices $i,j\in\{1,2\}$.
Given $\psi\in\Dom(Q)$, we have
$$
  Q[\psi]
  = \|h^{-1}\partial_1\psi\|_\Hilbert^2
  + \|\partial_2\psi\|_\Hilbert^2
  \,.
$$
Under the stated assumptions,
it is clear that the form domain of~$-\Delta_D^\Omega$
is just the Sobolev space
$
  \Sobi\big(\Real\times(-a,a)\big)
$.
If~$\mathcal{L}$ is injective,
then~$-\Delta_D^\Omega$ is nothing else
than the Dirichlet Laplacian
defined in the open subset
$
  \mathcal{L}\big(\Real\times(-a,a)\big)
$
of the ruled surface~(\ref{surface})
and expressed in the ``coordinates''~$(s,t)$.

\section{Geodesic strips}\label{Sec.geodesic}
%
The ruled strip~$\Omega$ is called a \emph{geodesic strip}
and denoted by~$\Omega_0$ if the reference curve~$\Gamma$
is a geodesic on~$\mathcal{L}$.
Since $\kappa\cos\theta$ is the geodesic curvature of~$\Gamma$
(when the latter is considered as a curve on~$\mathcal{L}$),
it is clear that~$\Omega$ is a geodesic strip
provided~$\Gamma$ is either a straight line (\ie~geodesic in~$\Real^3$)
or the straight lines
$
  t \mapsto \mathcal{L}(s,t) - \Gamma(s)
$
generating the ruled surface~(\ref{surface})
are tangential to the binormal vector field for each fixed~$s$.
The metric~(\ref{metric}) corresponding to~$\Omega_0$
acquires the form
\begin{equation}\label{metric.g}
  G_0 :=
  \begin{pmatrix}
    h_0^2 & 0 \\
    0   & 1
  \end{pmatrix}
  \,, \qquad
  h_0(s,t)
  := \sqrt{1 + t^2 \;\! \big[\tau(s)-\dot{\theta}(s)\big]^2}
  \,,
\end{equation}
and we denote by~$\Hilbert_0$, $Q_0$ and $-\Delta_D^{\Omega_0}$,
respectively, the corresponding Hilbert space
defined in analogy to~(\ref{Hilbert}),
the corresponding quadratic form defined in analogy to~(\ref{form})
and the associated Dirichlet Laplacian in~$\Omega_0$.

If $\tau-\dot{\theta}$ is equal to zero identically,
\ie~$\Omega_0$ is a flat surface due to~(\ref{Gauss}),
it is easy to see that the spectrum of~$-\Delta_D^{\Omega_0}$
coincides with the interval $[E_1,\infty)$,
where
$$
  E_1 := \pi^2/(2a)^2
$$
is the lowest eigenvalue of the Dirichlet Laplacian in $(-a,a)$.
In this section, we prove that
the presence of Gauss curvature
leads to a Hardy inequality
for the difference $-\Delta_D^{\Omega_0}-E_1$,
which has important consequences
for the stability of spectrum.
\begin{Theorem}\label{Thm.Hardy}
Given a positive number~$a$
and bounded continuous functions~$\tau$ and~$\dot{\theta}$,
let~$\Omega_0$ be the Riemannian manifold
$
  \big(\Real\times(-a,a),G_0\big)
$
with the metric given by~(\ref{metric.g}).
Assume that $\tau-\dot{\theta}$ is not identically zero
and that $a\,\|\tau-\dot{\theta}\|_\infty < \sqrt{2}$.
Then, for all $\psi\in\Sobi\big(\Real\times(-a,a)\big)$
and any~$s_0$ such that $(\tau-\dot{\theta})(s_0)\not=0$,
we have
$$
  Q_0[\psi] - E_1 \;\! \|\psi\|_{\Hilbert_0}^2
  \ \geq \
  c \, \|\rho^{-1}\psi\|_{\Hilbert_0}^2
  \qquad\mbox{with}\qquad
  \rho(s,t) := \sqrt{1+(s-s_0)^2}
  \,.
$$
Here~$c$ is a positive constant which depends
on~$s_0$, $a$ and $\tau-\dot{\theta}$.
\end{Theorem}

It is possible to find an explicit lower bound
for the constant~$c$;
we give an estimate in~(\ref{constant}) below.

Theorem~\ref{Thm.Hardy} implies
that the presence of Gauss curvature represents
a repulsive interaction in the sense that there
is no spectrum below~$E_1$ for all small potential-type
perturbations having $\mathcal{O}(s^{-2})$ decay at infinity.
Moreover, in the following Section~\ref{Sec.curved},
we show that this is also the case
for appropriate perturbations of the metric~(\ref{metric.g}).

In order to prove Theorem~\ref{Thm.Hardy},
we introduce the function $\lambda:\Real\to\Real$ by:
\begin{equation}\label{lambda}
  \lambda(s)
  := \inf_{\varphi\in\Smooth_0^\infty((-a,a))\setminus\{0\}} \,
  \frac{\int_{-a}^a |\dot{\varphi}(t)|^2 \, h_0(s,t) \, dt}
  {\int_{-a}^a |\varphi(t)|^2 \, h_0(s,t) \, dt}
  \ - \ E_1
\end{equation}
and keep the same notation for the function $\lambda \otimes 1$
on $\Real\times(-a,a)$. We have
\begin{Lemma}\label{Lem.crucial}
Under the hypotheses of Theorem~\ref{Thm.Hardy},
$\lambda$ is a continuous non-negative function
which is not identically equal to zero.
\end{Lemma}
\begin{proof}
For any fix~$s\in\Real$, we make the change of test function
$\phi:=\sqrt{h_0(s,\cdot)}\,\varphi$,
integrate by parts and arrive at
$$
  \lambda(s)
  = \inf_{\phi\in\Smooth_0^\infty((-a,a))\setminus\{0\}}
  \frac{\int_{-a}^a
  \left(
  |\dot{\phi}(t)|^2 - E_1\;\!|\phi(t)|^2 + V(s,t)\;\!|\phi(t)|^2
  \right) dt}
  {\int_{-a}^a |\phi(t)|^2 \, dt}
$$
with
\begin{equation}\label{potential}
  V(s,t) :=
  \frac{[\tau(s)-\dot{\theta}(s)]^2
  \left(2-t^2\,[\tau(s)-\dot{\theta}(s)]^2\right)}
  {4\,h_0(s,t)^4}
  \,.
\end{equation}
Under the hypotheses of Theorem~\ref{Thm.Hardy},
the function~$V$ is clearly continuous, non-negative
and not identically zero.
These facts together with the Poincar\'e inequality
$
  \int_{-a}^a |\dot{\phi}|^2
  \geq E_1 \int_{-a}^a |\phi|^2
$
valid for any $\phi\in\Smooth_0^\infty((-a,a))$
yield the claims of the Lemma.
\end{proof}

Assuming that the conclusion of Lemma~\ref{Lem.crucial} holds
and using the definition~(\ref{lambda}),
we get the estimate
\begin{equation}\label{local}
  Q_0[\psi] - E_1 \;\! \|\psi\|_{\Hilbert_0}^2
  \ \geq \
  \big\|h_0^{-1}\partial_1\psi\big\|_{\Hilbert_0}^2
  + \big\|\lambda^{1/2}\psi\big\|_{\Hilbert_0}^2
\end{equation}
valid for any
$
  \psi\in\Smooth_0^\infty\big(\Real\times(-a,a)\big)
$.
Neglecting the first term
on the right hand side of~(\ref{local}),
the inequality is already a Hardy inequality.
However, for applications it is more convenient
to replace the Hardy weight~$\lambda$ in~(\ref{local})
by the positive function~$c\;\!\rho^{-2}$
of Theorem~\ref{Thm.Hardy}.
This is possible by employing the contribution
of the first term based on
\begin{Lemma}\label{Lem.kinetic}
For any $\psi\in\Smooth_0^\infty\big(\Real\times(-a,a)\big)$,
\begin{multline*}
  \big(1+a^2\;\!\|\tau-\dot{\theta}\|_\infty^2\big)^{-1/2} \,
  \big\|\rho^{-1}\psi\big\|_{\Hilbert_0}^2
  \\
  \leq \
  16 \, \big(1+a^2\;\!\|\tau-\dot{\theta}\|_\infty^2\big)^{1/2} \,
  \big\|h_0^{-1}\partial_1\psi\big\|_{\Hilbert_0}^2
  + \big(2+64/|I|^2\big) \;\!
  \big\|\chi_I\;\!\psi\big\|_{\Hilbert_0}^2
  \,,
\end{multline*}
where~$I$ is any bounded subinterval of~$\Real$,
$\chi_I$~denotes the characteristic function
of the set $I\times(-a,a)$
and~$\rho$ is the function of Theorem~\ref{Thm.Hardy}
with~$s_0$ being the centre of~$I$.
\end{Lemma}
\begin{proof}
The Lemma is based on the following version
of the one-dimensional Hardy inequality:
\begin{equation}\label{one}
  \int_\Real \frac{|u(x)|^2}{x^{2}} \, dx
  \ \leq \
  4 \int_\Real |\dot{u}(x)|^2 \, dx
\end{equation}
valid for all $u\in\sobi(\Real)$ with $u(0)=0$.
Put $b:=|I|/2$.
We define the function $f:\Real\to\Real$ by
$$
  f(s) :=
  \begin{cases}
    1 & \mbox{for } \ |s-s_0| \geq b \,,
    \\
    |s-s_0|/b & \mbox{for } \ |s-s_0| < b \,,
  \end{cases}
$$
and keep the same notation for the function
$f \otimes 1$ on $\Real\times(-a,a)$.
For any $\psi\in\Smooth_0^\infty\big(\Real\times(-a,a)\big)$,
let us write
$
  \psi = f\psi + (1-f) \psi
$.
Applying~(\ref{one}) to the function
$s\mapsto(f\psi)(s,t)$ with~$t$ fixed,
we arrive at
\begin{eqnarray*}
  \int \frac{|\psi|^2}{\rho^2}
  & \leq &
  2 \int \frac{|f\psi|^2}{\rho^2-1}
  + 2 \int \chi_I \, |(1-f)\psi|^2
  \\
  & \leq &
  16 \int |\partial_1 f|^2 |\psi|^2
  + 16 \int |f|^2 |\partial_1\psi|^2
  + 2 \int \chi_I \, |(1-f)\psi|^2
  \\
  & \leq &
  16 \int |\partial_1\psi|^2
  + \big(2+16/b^2\big) \int \chi_I \, |\psi|^2
  \,,
\end{eqnarray*}
where the integration sign indicates
the integration over $\Real\times(-a,a)$.
Recalling the definition of~$\Hilbert_0$
and using the estimates
$$
  1
  \ \leq \ h_0^2 \ \leq \
  1+a^2\;\!\|\tau-\dot{\theta}\|_\infty^2
  \,,
$$
the Lemma follows at once.
\end{proof}

Now we are in a position to prove Theorem~\ref{Thm.Hardy}.
\begin{proof}[Proof of Theorem~\ref{Thm.Hardy}]
It suffices to prove the Theorem for functions~$\psi$
from the dense subspace $\Smooth_0^\infty\big(\Real\times(-a,a)\big)$.
Assume the hypotheses of Theorem~\ref{Thm.Hardy}
so that the conclusion of Lemma~\ref{Lem.crucial} holds.
Let~$I$ be any closed interval
on which~$\lambda$ is positive.
Writing
$$
  \big\|\lambda^{1/2}\psi\big\|_{\Hilbert_0}^2
  = \epsilon \, \big\|\lambda^{1/2}\psi\big\|_{\Hilbert_0}^2
  + (1-\epsilon) \, \big\|\lambda^{1/2}\psi\big\|_{\Hilbert_0}^2
  \qquad \mbox{with} \qquad
  \epsilon \in (0,1] \,,
$$
neglecting the second term of this decomposition,
estimating the first one by an integral over $I\times(-a,a)$
and applying Lemma~\ref{Lem.kinetic},
the inequality~(\ref{local}) yields
\begin{eqnarray*}
  \lefteqn{
  Q_0[\psi] - E_1 \;\! \|\psi\|_{\Hilbert_0}^2 }
  &&
  \\
  && \geq \
  \left[1-16\,\epsilon\,\min_I\lambda\,\big(2+64/|I|^2\big)^{-1}
  \big(1+a^2\;\!\|\tau-\dot{\theta}\|_\infty^2\big)^{1/2}
  \right]
  \big\|h_0^{-1}\partial_1\psi\big\|_{\Hilbert_0}^2
  \\
  && \phantom{\geq} \ +
  \epsilon \, \min_I\lambda \,
  \big(2+64/|I|^2\big)^{-1}
  \big(1+a^2\;\!\|\tau-\dot{\theta}\|_\infty^2\big)^{-1/2} \,
  \big\|\rho^{-1}\psi\big\|_{\Hilbert_0}^2
  \,.
\end{eqnarray*}
Choosing~$\epsilon$ as the minimum between~$1$
and the value such that the first term on the right hand side
of the last estimate vanishes,
we get the claim of Theorem~\ref{Thm.Hardy} with
\begin{equation}\label{constant}
  c \,\geq\, \min\left\{
  \frac{\min_I\lambda}
  {\big(2+64/|I|^2\big)\,
  \big(1+a^2\;\!\|\tau-\dot{\theta}\|_\infty^2\big)^{1/2}} \,,
  \frac{1}{16\,\big(1+a^2\;\!\|\tau-\dot{\theta}\|_\infty^2\big)}
  \right\}
  .
\end{equation}
\end{proof}
%

\section{Mildly curved strips}\label{Sec.curved}
%
Recall that the spectrum of~$-\Delta_D^{\Omega_0}$
coincides with the interval $[E_1,\infty)$
provided the Gauss curvature~(\ref{Gauss})
vanishes everywhere in the geodesic strip~$\Omega_0$.
On the other hand, it was proved in~\cite{K1}
that~$-\Delta_D^\Omega$ always possesses a spectrum below~$E_1$
provided the Gauss curvature~(\ref{Gauss}) vanishes everywhere
but~$\Gamma$ is not a geodesic on~$\mathcal{L}$.
In this section, we use the Hardy inequality
of Theorem~\ref{Thm.Hardy} to show
that the presence of Gauss curvature prevents
the spectrum to descend even if~$\Gamma$ is mildly curved.
\begin{Theorem}\label{Thm.stability}
Given a positive number~$a$
and bounded continuous functions~$\kappa$,
$\tau$ and~$\dot{\theta}$,
let~$\Omega$ be the Riemannian manifold~(\ref{strip})
with the metric given by~(\ref{metric}).
Assume that $\tau-\dot{\theta}$ is not identically zero
and that $a\,\|\tau-\dot{\theta}\|_\infty < \sqrt{2}$.
Assume also that for all $s\in\Real$,
\begin{equation*}
  |\kappa(s) \cos\theta(s)|
  \ \leq \
  \eps(s) := \frac{\eps_0}{1+s^2}
  \qquad\mbox{with}\qquad
  \eps_0 \in [0,a^{-1})
  \,.
\end{equation*}
Then there exists a positive number~$C$
such that $\eps_0 \leq C$ implies
\begin{equation}\label{stability}
  -\Delta_D^{\Omega}
  \ \geq \ E_1
  \,.
\end{equation}
Here~$C$ depends on~$a$
and on the constants~$c$ and~$s_0$ of Theorem~\ref{Thm.Hardy}.
\end{Theorem}

As usual, the inequality~(\ref{stability})
is to be considered in the sense of forms.
Actually, a stronger, Hardy-type inequality holds true,
\cf~(\ref{Hardy.c}) below.

An explicit lower bound for the constant~$C$
is given by the estimates made in the proof
of Theorem~\ref{Thm.stability} below.

As a direct consequence of Theorem~\ref{Thm.stability},
we get that the spectrum $[E_1,\infty)$ is stable as a set
provided the difference $\tau-\dot{\theta}$
vanishes at infinity:
\begin{Corollary}\label{Corol.stability}
In addition to hypotheses of Theorem~\ref{Thm.stability},
assume that $\tau(s)-\dot{\theta}(s)$ tends to zero
as $|s|\to\infty$.
Then
$$
  \spec(-\Delta_D^{\Omega}) = [E_1,\infty)
  \,.
$$
\end{Corollary}
\begin{proof}
Following the proof of~\cite[Sec.~3.1]{ChDFK} or \cite[Sec.~5]{KKriz}
based on a general characterization of essential spectrum
adopted from~\cite{DDI}, it is possible to show that
the essential spectrum~$-\Delta_D^{\Omega}$
coincides with the interval $[E_1,\infty)$,
while Theorem~\ref{Thm.stability} ensures that there
is no spectrum below~$E_1$.
\end{proof}
\begin{proof}[Proof of Theorem~\ref{Thm.stability}]
Let~$\psi$ belong to $\Smooth_0^\infty\big(\Real\times(-a,a)\big)$.
The proof is based on
an algebraic comparison of
$Q[\psi]-E_1\;\!\|\psi\|_\Hilbert^2$
with $Q_0[\psi]-E_1\;\!\|\psi\|_{\Hilbert_0}^2$
and the usage of Theorem~\ref{Thm.Hardy}.
For every $(s,t)\in\Real\times(-a,a)$,
we have
$$
  f_-(s) :=
  \sqrt{1-\frac{a\eps(s)[2+a\eps(s)]}
  {1+a^2\|\tau-\dot{\theta}\|_\infty^2}}
  \leq \frac{h(s,t)}{h_0(s,t)} \leq
  \sqrt{1+a\eps(s)[2+a\eps(s)]}
  =: f_+(s)
  \,.
$$
Here the lower bound is well defined and positive
provided $\eps_0 \leq (3a)^{-1}$,
and both the bounds behave as
$1+\mathcal{O}(\eps(s))$ as $\eps_0 \to 0$;
we keep the same notation~$f_\pm$
for the functions $f_\pm \otimes 1$ on $\Real\times(-a,a)$.
Consequently,
\begin{eqnarray*}
  Q[\psi] - E_1\;\!\|\psi\|_\Hilbert^2
  & \geq &
  \int_{\Real\times(-a,a)}
  f_+^{-1} \, h_0^{-1} \, |\partial_1\psi|^2
  \\
  && + \int_\Real ds \, f_-(s)
  \int_{-a}^a dt \, h_0(s,t) \left(
  |\partial_2\psi(s,t)|^2 - E_1\;\! |\psi(s,t)|^2
  \right)
  \\
  && - E_1 \int_{\Real\times(-a,a)} (f_+-f_-) \, h_0 \, |\psi|^2
  \,.
\end{eqnarray*}
Since the term in the second line is non-negative
due to~(\ref{lambda}) and Lemma~\ref{Lem.kinetic},
we can further estimate as follows:
\begin{eqnarray*}
  Q[\psi] - E_1\;\!\|\psi\|_\Hilbert^2
  & \geq &
  \min\big\{f_+(0)^{-1},f_-(0)\big\} \left(
  Q_0[\psi] - E_1\;\!\|\psi\|_{\Hilbert_0}^2
  \right)
  \\
  && - E_1 \int_{\Real\times(-a,a)} (f_+-f_-) \, h_0 \, |\psi|^2
  \,.
\end{eqnarray*}
Using Theorem~\ref{Thm.Hardy},
we finally obtain
\begin{equation}\label{Hardy.c}
  Q[\psi] - E_1\;\!\|\psi\|_\Hilbert^2
  \ \geq \
  \big\|w^{1/2}\psi\big\|_{\Hilbert_0}^2
  \,,
\end{equation}
where
\begin{equation*}
  w(s,t) :=
  \frac{c \, \min\big\{f_+(0)^{-1},f_-(0)\big\}}{1+(s-s_0)^2}
  - E_1 \;\! \big[f_+(s)-f_-(s)\big]
\end{equation*}
is positive for all sufficiently small~$\eps_0$.
\end{proof}
%

\section*{Acknowledgement}
%
The work has been supported by the Czech Academy of Sciences
and its Grant Agency within the projects IRP AV0Z10480505 and A100480501.

%
{\small

\providecommand{\bysame}{\leavevmode\hbox to3em{\hrulefill}\thinspace}

}
\ \bigskip \\
Department of Theoretical Physics \\
Nuclear Physics Institute \\
Academy of Sciences \\
250\,68 \v{R}e\v{z} near Prague \\
Czech Republic
\smallskip \\
\emph{E-mail:} krejcirik@ujf.cas.cz
\end{document}